\documentclass[11pt,reqno]{amsart}
\usepackage{amsfonts,amssymb,amsmath,color}

\usepackage[pdftex,colorlinks,
citecolor=black,linkcolor=black]{hyperref}

\DeclareMathOperator{\compl}{\scriptscriptstyle\complement}

\newtheorem{theorem}{Theorem}

\newtheorem{lemma}[theorem]{Lemma}

\newtheorem{remark}[theorem]{Remark}
\newtheorem{example}[theorem]{Example}

\numberwithin{equation}{section} \numberwithin{theorem}{section}

\begin{document}

\title{The increasing families of sets\\ generated by self-dual clutters}

\author{Andrey O. Matveev}
\email{andrey.o.matveev@gmail.com}

\begin{abstract} Using the KKS inequalities, we establish bounds on the numbers of\/ $k$-sets in the increasing families generated by self-dual clutters (i.e., clutters $\mathcal{A}$ that coincide with the blockers~$\mathfrak{B}(\mathcal{A})$) on their ground set of even cardinality.
\end{abstract}

\maketitle

\pagestyle{myheadings}

\markboth{THE INCREASING FAMILIES GENERATED BY SELF-DUAL CLUTTERS}{A.O.~MATVEEV}

\thispagestyle{empty}

\setcounter{tocdepth}{3}

\section{Introduction}

We denote by $E_t$, where $t\geq 3$, the set of integers $[t]:=[1,\ldots,t]:=\{1,\ldots,t\}$. A finite collection of sets $\mathcal{F}$ is called a {\em family}. The family~$\mathbf{2}^{[t]}$ of {\em all\/} subsets of the set~$E_t$ is called the {\em power set\/} of~$E_t$. The {\em empty set\/} is denoted by~$\hat{0}$, and the {\em empty family\/} containing {\em no sets\/} is denoted by~$\emptyset$. We denote by $|\!\cdot\!|$ the cardinalities of sets, while the numbers of sets in families are denoted by~$\#\,\cdot$. Sometimes we say that the cardinality of a set~$A$ is the {\em size\/} of~$A$.

If $\mathcal{F}$ is a set family such that $\emptyset\neq\mathcal{F}\neq\{\hat{0}\}$, then the union~\mbox{$\mathrm{V}(\mathcal{F}):=\bigcup_{F\in\mathcal{F}}F$} is called the {\em vertex set\/} of~$\mathcal{F}$. {\em Any\/} finite set~$S$ such that~$S\supseteq\mathrm{V}(\mathcal{F})$, {\em fixed\/} for someone's research purposes, is called the~{\em ground set\/} of the family~$\mathcal{F}$.

Given a nonempty family~$\mathcal{F}\subseteq\mathbf{2}^{[t]}$ on its ground set~$E_t$, we let~$\mathcal{F}^{\,\compl}:=\{F^{\compl}\colon$ $F\in\mathcal{F}\}$ denote the family of their {\em complements}, where $F^{\compl}:=E_t-F$. We also associate with the family $\mathcal{F}$ the family
\begin{equation*}
\mathcal{F}^{\hspace{0.2mm}\ast}:=\{G^{\compl}\colon G\in \mathbf{2}^{[t]}-\mathcal{F}\}\; .
\end{equation*}

A family of sets $\mathcal{A}$, such that $\emptyset\neq\mathcal{A}\neq\{\hat{0}\}$, is called a {\em nontrivial clutter\/}
if for any indices~$i\neq j$ of members of the family~$\{A_1,\ldots,A_{\alpha}\}=:\mathcal{A}$ we have~$A_i\not\subseteq A_j$.

For a subset $A\subseteq E_t$, the {\em principal increasing family\/} of sets $\{A\}^{\triangledown}\subseteq\mathbf{2}^{[t]}$, {\em generated\/} by the {\em one-member\/} clutter $\{A\}$ on its {\em ground set\/} $E_t$, is defined by $\{A\}^{\triangledown}:=\{B\subseteq E_t\colon$ $B\supseteq A\}$. In particular, we have~$\{\hat{0}\}^{\triangledown}=\mathbf{2}^{[t]}$, and~$\{E_t\}^{\triangledown}=\{E_t\}$.

Given a nonempty clutter~$\mathcal{A}:=\{A_1,\ldots,A_{\alpha}\}\subset\mathbf{2}^{[t]}$, the {\em increasing family\/} of sets~$\mathcal{A}^{\triangledown}\subseteq\mathbf{2}^{[t]}$, {\em generated\/} by $\mathcal{A}$ on its {\em ground set\/}~$E_t$, is defined as the union~$\mathcal{A}^{\triangledown}:=\bigcup_{A\in\mathcal{A}}\{A\}^{\triangledown}$ of the principal increasing families~$\{A\}^{\triangledown}$.

A subset $B\subseteq E_t$ is a {\em blocking set\/} of a nontrivial clutter~$\mathcal{A}\subset\mathbf{2}^{[t]}$ if it holds $|B\cap A|>0$, for each member $A\in\mathcal{A}$. The {\em blocker\/} $\mathfrak{B}(\mathcal{A})$ of the clutter~$\mathcal{A}$ is defined to be the family of all {\em inclusion-minimal blocking sets\/} of $\mathcal{A}$; see, e.g., the monographs~\cite{Berge89,Bretto,CCZ,Co,CH,CFKLMPPS,FLSZ,G,GJ,GLS,HHH,HY-Total,HY,Jukna-E2,M-PROM,Mirsky,NI,NW,SchU,Sch-C,V}. Additional references can be found in~\cite[Pt~2]{M-SC-V}. The increasing family~$\mathfrak{B}(\mathcal{A})^{\triangledown}$ is by definition the family of {\em all\/} blocking sets of the clutter~$\mathcal{A}$.

Recall that for a nontrivial clutter $\mathcal{A}$ on its ground set~$E_t$ we have
\begin{equation*}
\mathfrak{B}(\mathcal{A})^{\triangledown}=(\mathcal{A}^{\triangledown})^{\ast}\; .
\end{equation*}

Clutters $\mathcal{A}$ with the property
\begin{align*}
\mathfrak{B}(\mathcal{A})&=\mathcal{A}\; ,\\
\intertext{or, equivalently,}
\mathfrak{B}(\mathcal{A})^{\triangledown}&=\mathcal{A}^{\triangledown}\; ,
\end{align*}
are called {\em self-dual} or {\em identically self-blocking}; see, e.g.,~\cite{ACL-2019,AP-2018}\cite[\S{}2.1]{Berge89}\cite[Ch.~9]{Jukna-E2}\cite[\S{}5.7]{M-PROM} on such clutters.

A clutter~$\mathcal{A}$ on its ground set~$E_t$ is self-dual if and only if we have
\begin{equation}
\label{eq:30}
(\mathcal{A}^{\triangledown})^{\ast}=\mathcal{A}^{\triangledown}\; ,
\end{equation}
and although the (lack of) self-duality of clutters does not depend structurally on the cardinalities of their vertex sets and ground sets, relation~(\ref{eq:30}) suggests the following criterion:

\begin{remark} {\rm(see~\cite[Cor. 5.28(i)]{M-PROM})}
A clutter\/ $\mathcal{A}\subset\mathbf{2}^{[t]}$ on its ground set~$E_t$ is~{\em self-dual} if and only if
\begin{equation}
\label{eq:18}
\#\mathcal{A}^{\triangledown}=2^{t-1}\; .
\end{equation}
\end{remark}
In Theorem~\ref{th:3} we consider condition~(\ref{eq:18}) from a perspective of constraints that are satisfied by the numbers of $k$-sets in the increasing families~$\mathcal{A}^{\triangledown}$ generated by self-dual clutters~$\mathcal{A}\subset\mathbf{2}^{[t]}$ on their ground set~$E_t$ of even cardinality $t$ .

\section{Long $f$- and $h$-vectors of set families}

Let us associate with each family $\mathcal{F}\subseteq\mathbf{2}^{[t]}$ on the ground set~$E_t$ its \mbox{\em long\/ $f$-vector}
\begin{equation*}
\boldsymbol{f}(\mathcal{F};t):=\bigl(f_0(\mathcal{F};t),f_1(\mathcal{F};t),\ldots,f_t(\mathcal{F};t)\bigr)\in\mathbb{N}^{t+1}\; ,
\end{equation*}
where
\begin{equation*}
f_k(\mathcal{F};t):=\#\{F\in\mathcal{F}\colon |F|=k\}\; ,\ \ \ 0\leq k\leq t\; .
\end{equation*}
If `$\mathrm{x}$' is a formal variable, then the \mbox{\em long\/ $h$-vector}
\begin{equation*}
\boldsymbol{h}(\mathcal{F};t):=\bigl(h_0(\mathcal{F};t),h_1(\mathcal{F};t),\ldots,h_t(\mathcal{F};t)\bigr)\in\mathbb{Z}^{t+1}
\end{equation*}
of the family $\mathcal{F}$ is defined by means of the relation
\begin{equation*}
\sum_{i=0}^t h_i(\mathcal{F};t)\cdot\mathrm{x}^{t-i}:=\sum_{i=0}^t f_i(\mathcal{F};t)\cdot(\mathrm{x}-1)^{t-i}\; .
\end{equation*}

\begin{example} \label{th:4} {\rm(}cf.~Example~{\rm\ref{th:5}}{\rm)}
For any element\/ $a\in E_t$, the principal increasing family $\{\{a\}\}^{\triangledown}\subset\mathbf{2}^{[t]}$, generated by the {\em self-dual\/} clutter~$\{\{a\}\}$ on its ground set\/ $E_t$, is described by the vectors
\begin{align*}
\boldsymbol{f}(\{\{a\}\}^{\triangledown};t)&=(0,\tbinom{t-1}{1-1},\tbinom{t-1}{2-1},\ldots,\tbinom{t-1}{t-1})\\
\intertext{and}
\boldsymbol{h}(\{\{a\}\}^{\triangledown};t)&=(0,\phantom{-}1,\phantom{-}0,\ldots,0)\; .
\end{align*}
\end{example}

\begin{example}
{\rm(i)} {\rm(a)} Given the self-dual clutter
\begin{equation}
\label{eq:5}
\mathcal{A}:=\{\{1,2\},\{1,3\},\{2,3\}\}
\end{equation}
that generates on its ground set\/ $E_3=\mathrm{V}(\mathcal{A})$ the increasing family
\begin{equation*}
\mathcal{A}^{\triangledown}=\{\{1,2\},\{1,3\},\{2,3\},E_3\}\; ,
\end{equation*}
we have
\begin{align*}
\boldsymbol{f}(\mathcal{A}^{\triangledown};3)=(0,\phantom{-}0,\phantom{-}3,\phantom{-}1)\; ,\\
\boldsymbol{h}(\mathcal{A}^{\triangledown};3)=(0,\phantom{-}0,\phantom{-}3,-2)\; .
\end{align*}

{\rm(b)} If we choose the set $E_4\supsetneqq\mathrm{V}(\mathcal{A})$ as the ground set of the self-dual clutter~{\rm(\ref{eq:5})}, instead of the set~$E_3$, then $\mathcal{A}$ generates on $E_4$ the increasing family
\begin{equation*}
\mathcal{A}^{\triangledown}=\{\{1,2\},\{1,3\},\{2,3\},\{1,2,3\},\{1,2,4\},\{1,3,4\},\{2,3,4\},E_4\}\; ,
\end{equation*}
and we have
\begin{align*}
\boldsymbol{f}(\mathcal{A}^{\triangledown};4)=(0,\phantom{-}0,\phantom{-}3,\phantom{-}4,\phantom{-}1)\; ,\\
\boldsymbol{h}(\mathcal{A}^{\triangledown};4)=(0,\phantom{-}0,\phantom{-}3,-2,\phantom{-}0)\; .
\end{align*}

{\rm(ii)} {\rm(a)}
For the self-dual clutter
\begin{equation}
\label{eq:6}
\mathcal{A}:=\{\{2\}\}
\end{equation}
that generates on its ground set~$E_3\supsetneqq\mathrm{V}(\mathcal{A})$ the principal increasing family
\begin{equation*}
\mathcal{A}^{\triangledown}=\{\{2\},\{1,2\},\{2,3\},E_3\}\; ,
\end{equation*}
we have
\begin{align*}
\boldsymbol{f}(\mathcal{A}^{\triangledown};3)=(0,\phantom{-}1,\phantom{-}2,\phantom{-}1)\; ,\\
\boldsymbol{h}(\mathcal{A}^{\triangledown};3)=(0,\phantom{-}1,\phantom{-}0,\phantom{-}0)\; .
\end{align*}

{\rm(b)} For the self-dual clutter~{\rm(\ref{eq:6})} that generates on its ground set\newline \mbox{$E_4\supsetneqq\mathrm{V}(\mathcal{A})$} the principal increasing family
\begin{equation*}
\mathcal{A}^{\triangledown}=\{\{2\},\{1,2\},\{2,3\},\{2,4\},\{1,2,3\},\{1,2,4\},\{2,3,4\},E_4\}\; ,
\end{equation*}
we have
\begin{align*}
\boldsymbol{f}(\mathcal{A}^{\triangledown};4)=(0,\phantom{-}1,\phantom{-}3,\phantom{-}3,\phantom{-}1)\; ,\\
\boldsymbol{h}(\mathcal{A}^{\triangledown};4)=(0,\phantom{-}1,\phantom{-}0,\phantom{-}0,\phantom{-}0)\; .
\end{align*}

{\rm(iii)} The self-dual clutter
\begin{equation*}
\mathcal{A}:=\{\{1,2,3,4\},\{1,5\},\{2,5\},\{3,5\},\{4,5\}\}\; ,
\end{equation*}
that generates on its ground set~$E_5=\mathrm{V}(\mathcal{A})$ the increasing family
\begin{multline*}
\mathcal{A}^{\triangledown}=\{\{1,5\},\{2,5\},\{3,5\},\{4,5\},\{1,2,5\},\{1,3,5\},\{1,4,5\},\{2,3,5\},\\
\{2,4,5\},\{3,4,5\},\{1,2,3,4\},\{1,2,3,5\},\{1,2,4,5\},\\
\{1,3,4,5\},\{2,3,4,5\},E_5\}\; ,
\end{multline*}
is described by the vectors
\begin{align*}
\boldsymbol{f}(\mathcal{A}^{\triangledown};5)=(0,\phantom{-}0,\phantom{-}4,\phantom{-}6,\phantom{-}5,\phantom{-}1)\; ,\\
\boldsymbol{h}(\mathcal{A}^{\triangledown};5)=(0,\phantom{-}0,\phantom{-}4,-6,\phantom{-}5,-2)\; .
\end{align*}
\end{example}

\begin{remark} Given a family~$\mathcal{F}\subseteq\mathbf{2}^{[t]}$ on its ground set~$E_t$, we have{\rm:}
\begin{itemize}
\item[\rm(i)] {\rm(see~\cite[Prop.~2.1(iii)(a)]{M-PROM})}
\begin{align}
\label{eq:20}
h_{\ell}(\mathcal{F};t)&=(-1)^{\ell}\sum_{k=0}^{\ell}(-1)^k\tbinom{t-k}{t-\ell}\, f_k(\mathcal{F};t)\; , & 0\leq\;&\ell\leq t\; ;
\\
\label{eq:8}
f_{\ell}(\mathcal{F};t)&=\sum_{k=0}^{\ell}\tbinom{t-k}{t-\ell}\, h_k(\mathcal{F};t)\; , & 0\leq\;&\ell\leq t\; .
\end{align}

\item[\rm(ii)] {\rm(see~\cite[Prop.~2.1(iii)(b)]{M-PROM})}
\begin{align}
\nonumber
h_0(\mathcal{F};t)&=f_0(\mathcal{F};t)\; ;\\
\nonumber
h_1(\mathcal{F};t)&=f_1(\mathcal{F};t)-t f_0(\mathcal{F};t)\; ;\\
\nonumber
h_{t-1}(\mathcal{F};t)&=(-1)^{t-1}\sum_{k=0}^{t-1} (-1)^k (t-k)f_k(\mathcal{F};t)\; ;\\
\label{eq:22}
h_t(\mathcal{F};t)&=(-1)^t\sum_{k=0}^t (-1)^k f_k(\mathcal{F};t)\; ;\\
\nonumber
\sum_{k=0}^t h_k(\mathcal{F};t)&=f_t(\mathcal{F};t)\; .
\end{align}

\item[\rm(iii)] {\rm(see~\cite[Prop.~2.1(iii)(c)]{M-PROM})}
\begin{equation*}
\sum_{k=0}^t f_k(\mathcal{F};t)=\sum_{k=0}^t 2^{t-k}\, h_k(\mathcal{F};t)=\#\mathcal{F}\; .
\end{equation*}

\item[\rm(iv)]
\begin{equation*}
\boldsymbol{h}(\mathcal{F};t)+\boldsymbol{h}(\mathbf{2}^{[t]}-\mathcal{F};t)=(1,0,0,\ldots,0)=\boldsymbol{h}(\mathbf{2}^{[t]};t)\; .
\end{equation*}
\end{itemize}
\end{remark}

Being {\em redundant\/} analogues of {\em standard\/} $h$-vectors of abstract simplicial complexes~\cite{B-H,H-H,Mc,P-EN,St-96,Z}, {\em long\/} $h$-vectors (see, e.g.,~\cite[p.~170]{Mc-Sh}\cite[p.~265]{Mc-W}) demonstrate their usefulness below in relations~(\ref{eq:19}), where  they describe a curious enumerative connection between the families~$\mathcal{F}$ and~$\mathcal{F}^{\hspace{0.2mm}\ast}$.
\begin{remark}  Given a family~$\mathcal{F}\subseteq\mathbf{2}^{[t]}$ on its ground set~$E_t$,
we have{\rm:}
\begin{itemize}
\item[\rm(i)]
\begin{align*}
\#\mathcal{F}^{\hspace{0.2mm}\ast}+\#\mathcal{F}&=2^t\; .\\
\intertext{More precisely,}
f_{\ell}(\mathcal{F}^{\hspace{0.2mm}\ast};t)+f_{t-\ell}(\mathcal{F};t)&=\tbinom{t}{\ell}\; ,\ \ \ 0\leq \ell\leq t\; .
\end{align*}

\item[\rm(ii)] {\rm(see~{\rm\cite[Prop.~2.1(iii)(d)]{M-PROM}})}
\begin{equation}
\label{eq:19}
h_{\ell}(\mathcal{F}^{\hspace{0.2mm}\ast};t)+(-1)^{\ell}\sum_{k=\ell}^t\tbinom{k}{\ell}\,h_k(\mathcal{F};t)=\delta_{\ell,0}\; ,\ \ \
0\leq\ell\leq t\; ,
\end{equation}
where $\delta_{\ell,0}:=1$ if\/ $\ell=0$, and\/ $\delta_{\ell,0}:=0$ otherwise. Note that
\begin{equation*}
h_t(\mathcal{F}^{\hspace{0.2mm}\ast};t)=(-1)^{t+1}h_t(\mathcal{F};t)\; .
\end{equation*}
\end{itemize}
\end{remark}

\section{Abstract simplicial complexes $\Delta$, such that $\Delta^{\!\ast}=\Delta$, on their ground sets~$E_t$ of even cardinality $t$}

Although in the present paper we are concerned with enumerative properties of specific increasing families of sets, in this section we turn to
specific {\em abstract simplicial complexes\/} (i.e., ``{\em decreasing}'' families of ``{\em{}faces}'') because the numbers of faces of size~$k$ in complexes are known to satisfy the classical {\em Kruskal--Katona--Sch\"{u}tzenberger\/} (\!{\em{}KKS}) {\em inequalities}; see, e.g.,~\cite[\S{}10.3]{P-EN} and~\cite{An,Berge89,Bo,B-H,F-T,G-P,G-M,H-H,Jukna-E2,Knuth-4A,St-96,V,W,Z}.

Let $\Delta\subseteq\mathbf{2}^{[t]}$, where $\emptyset\neq\Delta\neq\{\hat{0}\}$, be an abstract simplicial complex on its ground set~$E_t$, with its {\em vertex set\/} $\mathrm{V}(\Delta)\subseteq E_t$. By definition of complex, the following implications hold: $(B\in\Delta$, $A\subset B)$ $\Longrightarrow$ $A\in\Delta$. The {\em inclusion-maximal\/} faces of $\Delta$ are called its {\em facets}. Sometimes one says that a complex with {\em one\/} facet (i.e., the power set of a set) is a {\em simplex}. The \mbox{long~$f$-vector} of the complex~$\Delta$ has the form
\begin{equation*}
\boldsymbol{f}(\Delta;t):=\bigl(1,\underbrace{f_1(\Delta;t)}_{:=|\mathrm{V}(\Delta)|> 0},\ldots,\underbrace{f_{d(\Delta)}(\Delta;t)}_{> 0},0,\ldots,0\bigr)\; ,
\end{equation*}
where $d(\Delta):=\max\{k\in [t]\colon f_k(\Delta;t)>0\}$.

\begin{example} \label{th:5} {\rm(}cf.~Example~{\rm\ref{th:4}}{\rm)}
For any element~$a\in E_t$, the {\em simplex}~$\Delta:=\{F\colon F\subseteq(E_t-\{a\})\}=\Delta^{\ast}$, whose facet is the subset $(E_t-\{a\})\subset E_t$ of size~$(t-1)$, is described by the vectors
\begin{align*}
\boldsymbol{f}(\Delta;t)&=(\tbinom{t-1}{0},\ldots,\tbinom{t-1}{t-2},\tbinom{t-1}{t-1},0)\\
\intertext{and}
\boldsymbol{h}(\Delta;t)&=(1,-1,\phantom{-}0,\ldots,0)\; .
\end{align*}
\end{example}

\begin{remark} \label{th:6}
An abstract simplicial complex~$\Delta$ with {\em vertex set\/}~$\mathrm{V}(\Delta)$, such that
$\Delta$ {\em coincides\/} with its {\em Alexander dual complex}\/ $\Delta^{\!\!\vee}$, defined by
\begin{equation*}
\Delta^{\!\!\vee}:=\{\mathrm{V}(\Delta)-F\colon F\in\mathbf{2}^{\mathrm{V(\Delta)}}-\Delta\}
\end{equation*}
when~$\mathrm{V}(\Delta^{\!\!\vee})=\mathrm{V}(\Delta)$, is known as an {\em Alexander self-dual complex}. Note that
\begin{equation*}
\Delta=\Delta^{\!\!\vee}\ \ \ \Longleftrightarrow\ \ \ \#\Delta=2^{|\mathrm{V}(\Delta)|-1}\; .
\end{equation*}
\end{remark}
\noindent{}See, e.g.,~\cite{V} and~\cite{BT-2009} on {\em combinatorial Alexander duality}.

\begin{remark}
Suppose that $\Delta\subset\mathbf{2}^{[t]}$ is an abstract simplicial complex, such that
$\Delta^{\!\ast}=\Delta$, on its ground set~$E_t$. Then we have{\rm:}
\begin{itemize}
\item[\rm(i)] $\#\Delta=\sum_{k=0}^{t}f_k(\Delta;t)=\sum_{k=0}^{d(\Delta)}f_k(\Delta;t)=\sum_{k=0}^t 2^{t-k}\,h_k(\Delta;t)$ \mbox{$=2^{t-1}$.}

\item[\rm(ii)]
\begin{equation}
\label{eq:17}
f_{\ell}(\Delta;t)+f_{t-\ell}(\Delta;t)=\tbinom{t}{\ell}\; ,\ \ \  0\leq \ell\leq t\; .
\end{equation}

\noindent{}In particular, if\/ $t$ is {\em even}, then
\begin{equation*}
f_{t/2}(\Delta;t)=\tfrac{1}{2}\tbinom{t}{t/2}\; .
\end{equation*}
\end{itemize}
\end{remark}

\begin{example}
{\rm(i)} The abstract simplicial complex
\begin{equation*}
\Delta:=\{\hat{0},\{1\},\{2\},\{3\},\{4\},\{1,2\},\{1,3\},\{2,3\}\}\; ,
\end{equation*}
such that $\Delta^{\!\ast}=\Delta$, on its ground set~$E_4=\mathrm{V}(\Delta)$, is described by the vectors
\begin{align*}
\boldsymbol{f}(\Delta;4)=(1,\phantom{-}4,\phantom{-}3,\phantom{-}0,\phantom{-}0)\; ,\\
\boldsymbol{h}(\Delta;4)=(1,\phantom{-}0,-3,\phantom{-}2,\phantom{-}0)\; .
\end{align*}
The complex~$\Delta$ is an Alexander self-dual complex, that is, $\Delta=\Delta^{\!\!\vee}$,
because~$8=\#\Delta=2^{|\mathrm{V}(\Delta)|-1}=2^{4-1}$.

{\rm(ii)} The simplex
\begin{equation*}
\Delta:=\{\hat{0},\{2\},\{3\},\{4\},\{2,3\},\{2,4\},\{3,4\},\{2,3,4\}\}\; ,
\end{equation*}
such that $\Delta^{\!\ast}=\Delta$, on its ground set~$E_4\supsetneqq\mathrm{V}(\Delta)$, is described by the vectors
\begin{align*}
\boldsymbol{f}(\Delta;4)=(1,\phantom{-}3,\phantom{-}3,\phantom{-}1,\phantom{-}0)\; ,\\
\boldsymbol{h}(\Delta;4)=(1,-1,\phantom{-}0,\phantom{-}0,\phantom{-}0)\; .
\end{align*}

\end{example}

\noindent$\bullet$ Suppose again that $\Delta\subset\mathbf{2}^{[t]}$ is an abstract simplicial complex, such that~\mbox{$\Delta^{\!\ast}=\Delta$,} on its ground set~$E_t$ of {\em even\/} cardinality~$t$. We would like to apply the KKS inequalities to
the complex~$\Delta$.

$\circ$ Since $\tfrac{1}{2}\tbinom{t}{t/2}=\tbinom{t-1}{t/2}$, we have the simplest $(t/2)$-{\em{}binomial expansion}
\begin{equation*}
f_{t/2}(\Delta;t)=\tbinom{\alpha_{t/2}(t/2)}{t/2}
\end{equation*}
of the middle component of the long $f$-vector $\boldsymbol{f}(\Delta;t)$, where
\begin{equation*}
\alpha_{t/2}(t/2)=t-1\; .
\end{equation*}
Then relations~(\ref{eq:17}) and the~KKS inequalities (lower and upper bound versions)~\cite[Cor.~10.1, Cor.~10.2]{P-EN} imply that
\begin{align*}
\tbinom{t}{(t/2)-1}-f_{(t/2)+1}(\Delta;t)=f_{(t/2)-1}(\Delta;t)&\geq\tbinom{\alpha_{t/2}(t/2)}{(t/2)-1}=\tbinom{t-1}{(t/2)-1}=\tbinom{t-1}{t/2}\; ,\\
\tbinom{t}{(t/2)+1}-f_{(t/2)-1}(\Delta;t)=f_{(t/2)+1}(\Delta;t)&\leq\tbinom{\alpha_{t/2}(t/2)}{(t/2)+1}=\tbinom{t-1}{(t/2)+1}=\tbinom{t-1}{(t/2)-2}\; .
\end{align*}
Further, note that if our complex $\Delta$ had the minimum possible number of its faces of size~$((t/2)-1)$, namely $\tbinom{t-1}{(t/2)-1}$ faces, then the following implication would hold:
\begin{multline*}
f_{(t/2)-1}(\Delta;t):=\tbinom{t-1}{(t/2)-1}=\tbinom{t-1}{t/2}\\
\Longrightarrow\ \ \
f_{(t/2)+1}(\Delta;t)=\tbinom{t}{(t/2)+1}-\tbinom{t-1}{(t/2)-1}\\
=\tbinom{t}{(t/2)-1}-\tbinom{t-1}{(t/2)-1}=
\tbinom{t-1}{(t/2)-2}\; .
\end{multline*}
Note also that $\tbinom{t}{(t/2)-1}-\tbinom{t-1}{(t/2)-1}=\tbinom{t-1}{(t/2)-2}=\tbinom{t-1}{(t/2)+1}$.

$\circ$ If $((t/2)-1)>1$, then suppose that our complex $\Delta$ indeed has the {\em minimum possible\/}
number of faces of size~$((t/2)-1)$, that is, $f_{(t/2)-1}(\Delta;t)$ $:=\tbinom{t-1}{(t/2)-1}$. In this case we have the $((t/2)-1)$-{\em{}binomial expansion}
\begin{equation*}
f_{(t/2)-1}(\Delta;t)=\tbinom{\alpha_{(t/2)-1}((t/2)-1)}{(t/2)-1}
\end{equation*}
of the component~$f_{(t/2)-1}(\Delta;t)$ of the vector~$\boldsymbol{f}(\Delta;t)$, where
\begin{equation*}
\alpha_{(t/2)-1}((t/2)-1)=t-1\; .
\end{equation*}
Relations~(\ref{eq:17}) and the~KKS inequalities (lower bound version)~\cite[Cor.~10.1]{P-EN} imply that
\begin{equation*}
\tbinom{t}{(t/2)-2}-f_{(t/2)+2}(\Delta;t)=f_{(t/2)-2}(\Delta;t)\geq\tbinom{\alpha_{(t/2)-1}((t/2)-1)}{(t/2)-2}=\tbinom{t-1}{(t/2)-2}\; .
\end{equation*}
On the other hand, suppose that the complex $\Delta$ has the {\em maximum possible\/}
number of faces of size~$((t/2)+1)$, that is, $f_{(t/2)+1}(\Delta;t):=\tbinom{t-1}{(t/2)+1}$. We have the~$((t/2)+1)$-{\em binomial expansion}
\begin{equation*}
f_{(t/2)+1}(\Delta;t)=\tbinom{\alpha_{(t/2)+1}((t/2)+1)}{(t/2)+1}
\end{equation*}
of the component~$f_{(t/2)+1}(\Delta;t)$ of the vector~$\boldsymbol{f}(\Delta;t)$, where
\begin{equation*}
\alpha_{(t/2)+1}((t/2)+1)=t-1\; .
\end{equation*}
Now relations (\ref{eq:17}) and the~KKS inequalities (upper bound version)~\cite[Cor.~10.2]{P-EN} imply that
\begin{equation*}
\tbinom{t}{(t/2)+2}-f_{(t/2)-2}(\Delta;t)=f_{(t/2)+2}(\Delta;t)\leq\tbinom{\alpha_{(t/2)+1}((t/2)+1)}{(t/2)+2}=\tbinom{t-1}{(t/2)+2}\; .
\end{equation*}
Note that $\tbinom{t}{(t/2)-2}-\tbinom{t-1}{(t/2)-2}=\tbinom{t-1}{(t/2)-3}=\tbinom{t-1}{(t/2)+2}$.

$\circ$ Proceeding by induction, we arrive at the following counterpart of~Theorem~\ref{th:3}:

\begin{lemma}
\label{th:2} {\rm(}see also Example~{\rm\ref{th:5}}{\rm)}
Let\/ $\Delta\subset\mathbf{2}^{[t]}$ be an abstract simplicial complex on its ground set\/ $E_t$ of\/ {\em even\/} cardinality\/ $t$.
If
\begin{equation*}
\Delta^{\!\ast}=\Delta\; ,
\end{equation*}
then we have
\begin{gather*}
f_0(\Delta;t)=1\; ,\ \ \text{and}\ \ \ f_t(\Delta;t)=0\; ;\\
f_{t/2}(\Delta;t)=\tfrac{1}{2}\tbinom{t}{t/2}=\tbinom{t-1}{t/2}\; ;\\
f_{(t/2)-i}(\Delta;t)\geq \tbinom{t-1}{(t/2)-i}=\tbinom{t-1}{(t/2)+i-1}\; ,\ \ \ 1\leq i\leq (t/2)-1\; ;\\
f_{(t/2)+j}(\Delta;t)\leq \tbinom{t-1}{(t/2)+j}=\tbinom{t-1}{(t/2)-j-1}\; ,\ \ \ 1\leq j\leq (t/2)-1\; ;\\
f_{(t/2)-k}(\Delta;t)+f_{(t/2)+k}(\Delta;t)=\tbinom{t}{(t/2)-k}=\tbinom{t}{(t/2)+k}\; ,\ \ \ 1\leq k\leq (t/2)-1\; .
\end{gather*}
\end{lemma}

\section{The increasing families of sets generated by self-dual clutters on their ground sets $E_t$ of even cardinality $t$ }

Let $\mathcal{F}\subset\mathbf{2}^{[t]}$ be a family on its ground set~$E_t$. The implication
\begin{equation*}
\mathcal{F}^{\hspace{0.2mm}\ast}=\mathcal{F}\ \ \ \ \Longleftrightarrow\ \ \ \ (\mathbf{2}^{[t]}-\mathcal{F})^{\ast}=\mathbf{2}^{[t]}-\mathcal{F}\; ,
\end{equation*}
allows us to obtain from Lemma~\ref{th:2} the following result:

\begin{theorem}
\label{th:3} {\rm(}see also Example~{\rm\ref{th:4}}{\rm)}
Let\/ $\mathcal{A}\subset\mathbf{2}^{[t]}$ be a clutter on its ground set\/ $E_t$ of\/ {\em even\/} cardinality\/ $t$. If the clutter\/ $\mathcal{A}$ is {\em self-dual}, that is,
\begin{equation*}
\mathfrak{B}(\mathcal{A})=\mathcal{A}\; ,
\end{equation*}
then we have
\begin{gather*}
f_0(\mathcal{A}^{\triangledown};t)=0\; ,\ \ \text{and}\ \ \ f_t(\mathcal{A}^{\triangledown};t)=1\; ;\\
f_{t/2}(\mathcal{A}^{\triangledown};t)=\tfrac{1}{2}\tbinom{t}{t/2}=\tbinom{t-1}{t/2}\; ;\\
f_{(t/2)-i}(\mathcal{A}^{\triangledown};t)\leq \tbinom{t-1}{(t/2)-i-1}
=\tbinom{t-1}{(t/2)+i}\; ,\ \ \ 1\leq i\leq (t/2)-1\; ;\\
f_{(t/2)+j}(\mathcal{A}^{\triangledown};t)\geq \tbinom{t-1}{(t/2)+j-1}
=\tbinom{t-1}{(t/2)-j}\; ,\ \ \ 1\leq j\leq (t/2)-1\; ;\\
f_{(t/2)-k}(\mathcal{A}^{\triangledown};t)+f_{(t/2)+k}(\mathcal{A}^{\triangledown};t)=\tbinom{t}{(t/2)-k}=\tbinom{t}{(t/2)+k}\; ,\ \ \ 1\leq k\leq (t/2)-1\; .
\end{gather*}
\end{theorem}

\section{Appendix: Remarks on set families $\mathcal{F}$ such that $\mathcal{F}^{\hspace{0.2mm}\ast}=\mathcal{F}$}

\begin{remark}
Let $\mathcal{F}\subset\mathbf{2}^{[t]}$ be a family, such that $\mathcal{F}^{\hspace{0.2mm}\ast}=\mathcal{F}$, on its ground set~$E_t$.
\begin{itemize}
\item[\rm(i)] We have
\begin{equation*}
\#\mathcal{F}=2^{t-1}\; .
\end{equation*}
More precisely,
\begin{equation}
\label{eq:28}
f_{\ell}(\mathcal{F};t)+f_{t-\ell}(\mathcal{F};t)=\tbinom{t}{\ell}\; ,\ \ \ 0\leq\ell\leq t\; .
\end{equation}

\item[\rm(ii)] Relations~{\rm(\ref{eq:28})} and~{\rm(\ref{eq:8})} imply that
\begin{equation*}
\sum_{k=0}^{\ell}\tbinom{t-k}{t-\ell}\, h_k(\mathcal{F};t)+\sum_{j=0}^{t-\ell}\tbinom{t-j}{\ell}\, h_j(\mathcal{F};t)=\tbinom{t}{\ell}\; , \ \ \ 0\leq\ell\leq t\; .
\end{equation*}
\end{itemize}
\end{remark}

Given a family~$\mathcal{F}\subset\mathbf{2}^{[t]}$, such that $\mathcal{F}^{\hspace{0.2mm}\ast}=\mathcal{F}$, on its ground set~$E_t$, from~(\ref{eq:20}) we have for any~$\ell$, where~$0\leq\ell\leq t$:
\begin{equation*}
\begin{split}
h_{\ell}(\mathcal{F};t)&=(-1)^{\ell}\sum_{k=0}^{\ell}(-1)^k\tbinom{t-k}{t-\ell}\, f_k(\mathcal{F};t)\\
&=(-1)^{\ell}\sum_{k=0}^{\ell}(-1)^k\tbinom{t-k}{t-\ell}\, \Bigl(\tbinom{t}{k}-f_{t-k}(\mathcal{F};t)\Bigr)\; .
\end{split}
\end{equation*}

\begin{remark} Let~$\mathcal{F}\subset\mathbf{2}^{[t]}$ be a family, such that $\mathcal{F}^{\hspace{0.2mm}\ast}=\mathcal{F}$, on its ground set~$E_t$.
\begin{itemize}
\item[\rm(i)] We have
\begin{align}
\label{eq:21}
h_{\ell}(\mathcal{F};t)&=\delta_{\ell,0}
-(-1)^{t-\ell}\sum_{j=t-\ell}^t(-1)^j\tbinom{j}{t-\ell}\, f_j(\mathcal{F};t)\; , & 0\leq\; &\ell\leq t\; ;\\
\label{eq:14}
f_{\ell}(\mathcal{F};t)&=\tbinom{t}{\ell}-\sum_{j=0}^{t-\ell}\tbinom{t-j}{\ell}\, h_j(\mathcal{F};t)\; , & 0\leq\; &\ell\leq t\; .
\end{align}

\item[\rm(ii)] Relations {\rm(\ref{eq:21})} and {\rm(\ref{eq:20})} yield
\begin{multline*}
(-1)^{\ell}\sum_{k=0}^{\ell}(-1)^k\tbinom{t-k}{t-\ell}\, f_k(\mathcal{F};t)\\
+(-1)^{t-\ell}\sum_{j=t-\ell}^t(-1)^j\tbinom{j}{t-\ell}\, f_j(\mathcal{F};t)=\delta_{\ell,0}\; ,\ \ \  0\leq\ell\leq t\; .
\end{multline*}
\end{itemize}
\end{remark}

\begin{remark} Let\/ $\mathcal{F}\subset\mathbf{2}^{[t]}$ be a family, such that $\mathcal{F}^{\hspace{0.2mm}\ast}=\mathcal{F}$, on its ground set~$E_t$ of\/ {\em odd\/} cardinality~$t$. From~{\rm(\ref{eq:22})} we have
\begin{equation*}
\begin{split}
h_t(\mathcal{F};t)&=(-1)^t\sum_{k=0}^t(-1)^k f_k(\mathcal{F};t)\\
&=-\sum_{k=0}^{\lfloor t/2 \rfloor}(-1)^k \Bigl(f_k(\mathcal{F};t)-(\tbinom{t}{k}-f_k(\mathcal{F};t))\Bigr)\\
&=\tbinom{t-1}{(t-1)/2}-2\sum_{k=0}^{(t-1)/2}(-1)^k f_k(\mathcal{F};t)\; ,
\end{split}
\end{equation*}
and
\begin{equation*}
\begin{split}
h_t(\mathcal{F};t)&=(-1)^t\sum_{j=0}^t(-1)^j f_j(\mathcal{F};t)\\
&=-\sum_{j=\lceil t/2\rceil}^t (-1)^j \Bigl(f_j(\mathcal{F};t)-(\tbinom{t}{j}-f_j(\mathcal{F};t))\Bigr)\\
&=-\tbinom{t-1}{(t-1)/2}-2\sum_{j=(t+1)/2}^t (-1)^j f_j(\mathcal{F};t)\; .
\end{split}
\end{equation*}
As a consequence, we see that
\begin{equation*}
\sum_{k=0}^{(t-1)/2}(-1)^k f_k(\mathcal{F};t)
=\tbinom{t-1}{(t-1)/2}+\sum_{j=(t+1)/2}^t (-1)^j f_j(\mathcal{F};t)\; .
\end{equation*}
\end{remark}

\begin{remark}
Let $\mathcal{F}\subset\mathbf{2}^{[t]}$ be a family, such that $\mathcal{F}^{\hspace{0.2mm}\ast}=\mathcal{F}$, on its ground set~$E_t$.
\begin{itemize}
\item[\rm(i)] From~{\rm(\ref{eq:19})} we have
\begin{equation}
\label{eq:23}
h_{\ell}(\mathcal{F};t)=\delta_{\ell,0}+(-1)^{\ell+1}\sum_{k=\ell}^t\tbinom{k}{\ell}\, h_k(\mathcal{F};t)\; ,\ \ \ 0\leq \ell\leq t\; .
\end{equation}
If\/ $t$ is {\em even}, then we see that
\begin{equation}
\label{eq:24}
h_t(\mathcal{F};t)=0\; .
\end{equation}

\item[\rm(ii)]
For {\em even} indices\/ $\ell$, where $2\leq\ell\leq t$, relations~{\rm(\ref{eq:23})} imply that
\begin{equation}
\label{eq:25}
\sum_{k=\ell}^t\tbinom{k}{\ell-1}\, h_k(\mathcal{F};t)=0\; ,
\end{equation}
that is,
\begin{equation*}
\ell\,h_{\ell}(\mathcal{F};t)+\sum_{k=\ell+1}^t\tbinom{k}{\ell-1}\, h_k(\mathcal{F};t)=0\; ;
\end{equation*}
we also have
\begin{equation*}
2\, h_{\ell}(\mathcal{F};t)+
\sum_{k=\ell+1}^t\tbinom{k}{\ell}\, h_k(\mathcal{F};t)=0\; .
\end{equation*}

If\/ $t$ is {\em odd}, then we have
\begin{equation*}
(t-1)\, h_{t-1}(\mathcal{F};t)+\tbinom{t}{2}\,h_t(\mathcal{F};t)=0\; ,
\end{equation*}
that is,
\begin{equation*}
2\, h_{t-1}(\mathcal{F};t)+t\, h_t(\mathcal{F};t)=0\; .
\end{equation*}

If\/ $t$ is {\em even}, then we have
\begin{align}
\nonumber
\tfrac{t}{2}\, h_{t/2}(\mathcal{F};t)+
\sum_{k=(t/2)+1}^{t-1}\tbinom{k}{(t/2)-1}\, h_k(\mathcal{F};t)&=0\; ,\\
\intertext{and}
\label{eq:27}
h_{t/2}(\mathcal{F};t)+\tfrac{1}{2}
\sum_{k=(t/2)+1}^{t-1}\tbinom{k}{t/2}\, h_k(\mathcal{F};t)&=0\; .
\end{align}
\end{itemize}
\end{remark}

\begin{remark}
Let\/ $\mathcal{F}\subset\mathbf{2}^{[t]}$ be a family, such that $\mathcal{F}^{\hspace{0.2mm}\ast}=\mathcal{F}$, on its ground set~$E_t$.

If\/ $t$ is {\em even}, then relations~~{\rm(\ref{eq:24})} and~{\rm(\ref{eq:25})} in the case $\ell:=2$ yield
\begin{equation*}
h_t(\mathcal{F};t)=0\; ,\ \ \ \text{and}\ \ \ \ \sum_{k=2}^{t-1}k\, h_k(\mathcal{F};t)=0\; .
\end{equation*}

Similarly, if\/ $t$ is {\em odd}, then from~{\rm(\ref{eq:25})} we have
\begin{equation*}
\sum_{k=2}^t k\, h_k(\mathcal{F};t)=0\; .
\end{equation*}
\end{remark}

\begin{remark}
Let\/ $\mathcal{F}\subset\mathbf{2}^{[t]}$ be a family, such that\/ $\mathcal{F}^{\hspace{0.2mm}\ast}=\mathcal{F}$, on its ground set\/ $E_t$ of {\em even}
cardinality\/ $t$.
\begin{itemize}
\item[\rm(i)]
On the one hand, since
\begin{equation}
\label{eq:29}
f_{t/2}(\mathcal{F};t)=\tfrac{1}{2}\tbinom{t}{t/2}\; ,\\
\end{equation}
relations~{\rm(\ref{eq:8})} imply that
\begin{equation}
\label{eq:9}
\begin{split}
\sum_{k=0}^{t/2}\tbinom{t-k}{t/2}\, h_k(\mathcal{F};t)&=f_{t/2}(\mathcal{F};t)\\
=h_{t/2}(\mathcal{F};t)+\sum_{k=0}^{(t/2)-1}\tbinom{t-k}{t/2}\, h_k(\mathcal{F};t)&=\tfrac{1}{2}\tbinom{t}{t/2}\; .
\end{split}
\end{equation}

On the other hand, {\rm(\ref{eq:29})} and relations~{\rm(\ref{eq:14})} imply that
\begin{equation*}
\tbinom{t}{t/2}-\sum_{j=0}^{t/2}\tbinom{t-j}{t/2}\, h_{t/2}(\mathcal{F};t)=f_{t/2}(\mathcal{F};t)
=\tfrac{1}{2}\tbinom{t}{t/2}\; ,
\end{equation*}
that is,
\begin{equation*}
\sum_{j=0}^{t/2}\tbinom{t-j}{t/2}\, h_{t/2}(\mathcal{F};t)=\tfrac{1}{2}\tbinom{t}{t/2}\; .
\end{equation*}

\item[\rm(ii)]
Relations~{\rm(\ref{eq:27})} and~{\rm(\ref{eq:9})} imply that
\begin{equation*}
\sum_{k=0}^{(t/2)-1}\tbinom{t-k}{t/2}\, h_k(\mathcal{F};t)
-\tfrac{1}{2}\sum_{j=(t/2)+1}^{t-1}\tbinom{j}{t/2}\, h_j(\mathcal{F};t)
=\tfrac{1}{2}\tbinom{t}{t/2}\; ,
\end{equation*}
and
\begin{equation*}
h_{t/2}(\mathcal{F};t)=
\tfrac{1}{4}\tbinom{t}{t/2}-\tfrac{1}{2}\sum_{k=0}^{(t/2)-1}\tbinom{t-k}{t/2}\, h_k(\mathcal{F};t)
-\tfrac{1}{4}\sum_{j=(t/2)+1}^{t-1}\tbinom{j}{t/2}\, h_j(\mathcal{F};t)
\; .
\end{equation*}
\end{itemize}
\end{remark}

\vspace{5mm}

\begin{thebibliography}{99.}
\bibitem{ACL-2019}
{\em Abdi~A.}, {\em Cornu\'{e}jols~G.}, {\em Lee~D.} Identically self-blocking clutters. 20th International Conference Integer Programming and Combinatorial Optimization IPCO 2019 held in Ann~Arbor, MI, USA, May 22--24, 2019. {\em A.~Lodi\/} and {\em V.~Nagarajan\/} (eds.).
Lecture Notes in Computer Science, 11480. Cham: Springer, 2019, pp.~1--12.

\bibitem{AP-2018}
{\em Abdi~A.}, {\em Pashkovich~K.} Delta minors, delta free clutters, and entanglement. SIAM Journal on Discrete Mathematics, 2018, 32, no.~3, pp.~1750--1774.

\bibitem{An}
{\em Anderson~I.} Combinatorics of finite sets. Corrected reprint of the 1989 edition. Mineola, NY: Dover Publications, Inc., 2002.

\bibitem{Berge89}
{\em Berge~C.} Hypergraphs. Combinatorics of finite sets. Translated from the French. North-Holland Mathematical Library, 45. Amsterdam: North-Holland Publishing Co., 1989.

\bibitem{BT-2009}
{\em Bj\"{o}rner~A.}, {\em Tancer~M.} Note: Combinatorial Alexander duality---A short and elementary proof.
Discrete and Computational Geometry, 2009, 42, no.~4, pp.~586--593.

\bibitem{Bo}
{\em Bollob\'{a}s~B.} Combinatorics. Set systems, hypergraphs, families of vectors and combinatorial probability. Cambridge: Cambridge University Press, 1986.

\bibitem{Bretto}
{\em Bretto~A.} Hypergraph theory. An introduction. Mathematical Engineering. Cham: Springer, 2013.

\bibitem{B-H}
{\em Bruns W.}, {\em Herzog~J}. Cohen--Macaulay rings. Second edition. Cambridge Studies in Advanced
Mathematics, 39. Cambridge: Cambridge University Press, 1998.

\bibitem{CCZ}
{\em Conforti~M.}, {\em Cornu\'{e}jols~G.}, {\em Zambelli~G.} Integer programming. Cham: Springer, 2014.

\bibitem{Co}
{\em Cornu\'{e}jols~G.} Combinatorial optimization. Packing and covering. CBMS--NSF Regional
Conference Series in Applied Mathematics, 74. Philadelphia, PA: Society for Industrial and
Applied Mathematics (SIAM), 2001.

\bibitem{CH}
{\em Crama~Y.}, {\em Hammer~P.L.} Boolean functions. Theory, algorithms, and applications. Encyclopedia of Mathematics and its Applications, 142. Cambridge: Cambridge University Press, 2011.

\bibitem{CFKLMPPS}
{\em Cygan~M.}, {\em Fomin~F.V.}, {\em Kowalik~{\L}.}, {\em Lokshtanov~D.}, {\em Marx~D.}, {\em Pilipczuk~M.}, {\em Pilipczuk~M.}, {\em Saurabh~S.} Parameterized algorithms. Cham: Springer, 2015.

\bibitem{FLSZ}
{\em Fomin~F.V.}, {\em Lokshtanov~D.}, {\em Saurabh~S.}, {\em Zehavi~M.} Kernelization. Theory of parameterized preprocessing. Cambridge: Cambridge University Press, 2019.

\bibitem{F-T}
{\em Frankl~P.}, {\em Tokushige~N.} Extremal problems for finite sets. Student Mathematical Library, 86. Providence, RI: American Mathematical Society, 2018.

\bibitem{G}
{\em Gainanov~D.N.} Graphs for pattern recognition. Infeasible systems of linear inequalities. Berlin: De\,Gruyter, 2016.

\bibitem{GJ}
{\em Garey~M.R.}, {\em Johnson~D.S.} Computers and intractability. A guide to the theory of~NP-completeness. A Series of Books in the Mathematical Sciences. San~Francisco,~CA: W.H.~Freeman and~Co., 1979.

\bibitem{G-P}
{\em Gerbner~D,}, {\em Patk\'{o}s~B.} Extremal finite set theory. Discrete Mathematics and its Applications (Boca Raton). Boca Raton, FL: CRC Press, 2019.

\bibitem{G-M}
{\em Godsil~C.}, {\em Meagher~K.} Erd\H{o}s–Ko–Rado theorems: Algebraic approaches. Cambridge
Studies in Advanced Mathematics, 149. Cambridge: Cambridge University Press, 2016.

\bibitem{GLS}
{\em Gr\"{o}tschel~M.}, {\em Lov\'{a}sz~L.}, {\em Schrijver~A.} Geometric algorithms and combinatorial optimization.
Second edition. Algorithms and Combinatorics, 2. Berlin: Springer-Verlag, 1993.

\bibitem{HHH}
{\em Haynes~T.W.}, {\em Hedetniemi~S.T.}, {\em Henning~M.A.} (eds.). Structures of domination in graphs. Developments in Mathematics, 66. Cham: Springer, 2021.

\bibitem{HY-Total}
{\em Henning~M.A.}, {\em Yeo~A.} Total domination in graphs. Springer Monographs in Mathematics. New~York: Springer, 2013.

\bibitem{HY}
{\em Henning~M.A.}, {\em Yeo~A.} Transversals in linear uniform hypergraphs. Developments in Mathematics, 63. Cham: Springer, 2020.

\bibitem{H-H}
{\em Herzog~J.}, {\em Hibi~T.} Monomial ideals. Graduate Texts in Mathematics, 260. London: Springer-Verlag London, Ltd., 2011.

\bibitem{Jukna-E2}
{\em Jukna~S.} Extremal combinatorics. With applications in computer science. Second edition. Texts in Theoretical Computer Science. An EATCS Series. Heidelberg: Springer, 2011.

\bibitem{Knuth-4A}
{\em Knuth~D.E.} The art of computer programming. Vol.~4A. Combinatorial algorithms. Part~1. Upper~Saddle~River, NJ: Addison-Wesley, 2011.

\bibitem{M-PROM}
{\em Matveev~A.O.} Pattern recognition on oriented matroids. Berlin: De\,Gruyter, 2017.

\bibitem{M-SC-V}
{\em Matveev A.O.} Pattern recognition on oriented matroids: Symmetric cycles in the
hypercube graphs. V. Preprint [arXiv:2106.03832], 2021.

\bibitem{Mc}
{\em McMullen~P.} The numbers of faces of simplicial polytopes. Israel Journal of Mathematics, 1971, 9, 559--570.

\bibitem{Mc-Sh}
{\em McMullen~P.}, {\em Shephard G.C.} Convex polytopes and the upper bound conjecture. Prepared in
collaboration with {\em J.E.\,Reeve} and {\em A.A.\,Ball}. London Mathematical Society Lecture Note Series,
vol.~3. London: Cambridge University Press, 1971.

\bibitem{Mc-W}
{\em McMullen P.}, {\em Walkup~D.W.} A generalized lower-bound conjecture for simplicial polytopes.
Mathematika, 1971, 18, pp. 264--273.

\bibitem{Mirsky}
{\em Mirsky~L.} Transversal theory. An account of some aspects of combinatorial mathematics. Mathematics in Science and Engineering, 75. New~York: Academic Press, 1971.

\bibitem{NI}
{\em Nagamochi~H.}, {\em Ibaraki~T.} Algorithmic aspects of graph connectivity. Encyclopedia of Mathematics and its Applications, 123. Cambridge: Cambridge University Press, 2008.

\bibitem{NW}
{\em Nemhauser~G.L.}, {\em Wolsey~L.A.} Integer and combinatorial optimization. Reprint of the~1988 original. Wiley-Interscience Series in Discrete Mathematics and Optimization. New York: John Wiley \& Sons, Inc., 1999.

\bibitem{P-EN}
{\em Petersen T.K.} Eulerian numbers. With a foreword by {\em R.P.~Stanley}. Birkh\"{a}user Advanced Texts: Basler Lehrb\"{u}cher.  New York: Birkh\"{a}user/Springer, 2015.

\bibitem{SchU}
{\em Scheinerman~E.R.}, {\em Ullman~D.H.} Fractional graph theory. A rational approach to the theory of graphs. With a foreword by\! {\em C.~Berge}. Reprint of the 1997 original. Mineola, NY: Dover Publications, Inc., 2011.

\bibitem{Sch-C}
{\em Schrijver~A.} Combinatorial optimization. Polyhedra and efficiency. Vol.~C. Disjoint paths,
hypergraphs. Chapters 70--83. Algorithms and Combinatorics, 24,C. Berlin: Springer-Verlag,
2003.

\bibitem{St-96}
{\em Stanley R.P.} Combinatorics and commutative algebra. Second edition. Progress in Mathematics, 41. Boston, MA: Birkhäuser, 1996.

\bibitem{V}
{\em Villarreal~R.H.} Monomial algebras. Second edition. Monographs and Research Notes in Mathematics. Boca~Raton, FL: CRC Press, 2015.

\bibitem{W}
{\em West~D.B.} Combinatorial mathematics. Cambridge: Cambridge University Press, 2021.

\bibitem{Z}
{\em Ziegler~G.M.} Lectures on polytopes. Revised edition. Graduate Texts in Mathematics, 152. Berlin: Springer-Verlag, 1998.

\end{thebibliography}
\end{document}